\documentclass[11pt,a4paper,reqno]{amsart}
\usepackage{enumerate,array}
\usepackage{oppi}
\newtheorem{theorem}{Theorem}[section]
\newtheorem{lemma}[theorem]{Lemma}

\theoremstyle{definition}

\theoremstyle{remark}
\newtheorem{remark}[theorem]{Remark}

\numberwithin{equation}{section}
\renewcommand\Eq{\club\club\show \Eq}
\renewcommand\equ{\club\club\show \equ}
\renewcommand\thm{\club\club\show \thm}

\newcommand\lmax{{\l_{\rm max}}}

\newcommand\repres[1]{\mathbf{#1}}
\newcommand\rT{\repres{T}}
\newcommand\rR{\repres{R}}
\newcommand\rY{\repres{Y}}
\newcommand\rD{\repres{D}}
\newcommand\hZ{\widehat Z}
\newcommand\Smm{\operatorname{Sym}}
\newcommand\symCG{\Smm \bC G}
\newcommand\symCGp{(\Smm \bC G)_+}

\newcommand\Sc{Schreier}
\newcommand\secskip{\bigskip}
\begin{document} 
%

\title[Cayley graphs generated by initial reversals]{%
Cayley graphs on the symmetric group generated by initial reversals
have unit spectral gap
}
%
%
%
%
%
%
\author[Filippo Cesi]{Filippo Cesi}
\address{%
    Filippo Cesi\hfill\break
    \indent Dipartimento di Fisica\hfill\break
    \indent Universit\`a di Roma ``La Sapienza", Italy\hfil\break
    \indent and SMC, INFM-CNR.
}
\email{filippo.cesi@roma1.infn.it}

%
%
%
%
\begin{abstract}
In a recent paper Gunnells, Scott and Walden 
have determined the complete spectrum of the Schreier graph
on the symmetric group
corresponding to the Young subgroup $S_{n-2}\times S_2$
and generated by initial reversals.
In particular they find that the first nonzero eigenvalue, or spectral gap,
of the Laplacian is always 1, and report that ``empirical evidence''
suggests that this also holds for the corresponding Cayley graph.
We provide a simple proof of this last assertion, based on the
decomposition of the Laplacian of Cayley graphs,
into a direct sum of irreducible representation matrices
of the symmetric group.
\end{abstract}

\maketitle
\thispagestyle{empty}
%
%
%
%
%
\section{Introduction} 

\noindent
If $G$ is a finite group, 
$H$ is a subgroup of $G$ and $Z$ is a generating set of $G$,
we can construct the \Sc\ graph
$\cG=X(G,H,Z)$ whose vertices are the left-cosets $G/H$,
and whose edges are the pairs $(gH, zgH)$ with $gH\in G/H$ and $z\in Z$.
We assume that the generating set $Z$ is \textit{symmetric}, \ie $z\in Z$
if and only if $z^{-1}\in Z$. In this case the graph $X(G,H,Z)$ is undirected.
If $H=\{1\}$ we denote with $X(G,Z) = X(G,\{1\},Z)$ the
Cayley graph of $G$ associated to the generating set $Z$.
If $A_{\cG}$ is the \textit{adjacency matrix}
of $\cG$ and $\D_{\cG}$ the corresponding Laplacian,
since $\cG$ is $|Z|$--regular (counting loops), we have
\begin{equation*}
  \D_\cG = |Z| \, \id - A_\cG \,,
\end{equation*}
where $|Z|$ stands for the cardinality of the set $Z$.
The Laplacian is symmetric 
and positive-semidefinite, hence
its eigenvalues are real and nonnegative and can be ordered
as 
\begin{equation*}
0=\l_1(\D_\cG) \le \l_2(\D_\cG) \le\cdots\le \l_n(\D_\cG) \,.
\end{equation*}
Since $Z$ generates $G$, the graph $\cG$ is connected, which implies
that $0$ is a simple eigenvalue 
with constant eigenvector, while 
$\l_2(\D_\cG)$ is strictly positive.
The second eigenvalue of the Laplacian
is also called the \textit{spectral gap} of the graph $\cG$
and we denote it with $\gap \cG$.
For regular graphs, it coincides with the difference
between the two largest eigenvalues of the adjacency matrix.

In \cite{GSW} the authors consider the
\Sc\ graph $X(S_n, S_{(n-2,2)}, Z_n)$ where $S_n$ is the symmetric
group, $S_{(n-2,2)}$ is the Young subgroup corresponding
to the partition $(n-2,2)$, isomorphic to $S_{n-2}\times S_2$,
and $Z_n := \{ r_1, \ldots, r_n \}$, where
$r_k$ is the permutation which reverses the order of the
first $k$ positive integers
\begin{equation}\label{eq:rk}
  r_k : \{1,2,\ldots, n\} \too 
  \{ k,\, k-1,\, \ldots,\, 2,\,1,\,  k+1, \ldots, n\} \,.
\end{equation}
In \cite{GSW} the spectrum of the Laplacian was
determined, and it turns out that 
\begin{equation*}
  \spec \D_{X(S_n, S_{(n-2,2)}, Z_n)} \sset \{0,1,\ldots,n\}
\end{equation*}
\textit{with equality} if $n\ge 8$.
It was also proved that $\gap X(S_n, S_{(n-2,2)}, Z_n) = 1$
for all $n\ge 3$.
On the other side it was shown in \cite{Nas} that, if $Z_n$
is a set of reversals with $|Z_n| = o(n)$, then 
$\gap X(S_n, S_{(n-2,2)}, Z_n) \to 0$ as $n\to\oo$.
Hence, results in \cite{GSW} show that condition 
$|Z_n|=o(n)$ is, in a sense, optimal.

It is easy to see that if $K$ is a subgroup of $H$, then
the spectrum of $X(G,H,Z)$ is a subset of the spectrum of $X(G,K,Z)$.
In particular the spectrum of $X(G,H,Z)$ is a subset of the spectrum
of the Cayley graph $X(G,Z)$, thus we have
$\gap X(G,Z) \le \gap X(G,H,Z)$.
By consequence we get, for what concerns the symmetric group
generated by initial reversals,%
\footnote{when $n=2$ a trivial computation yields $\gap(S_2,Z_2)=2$.}
\begin{equation}\label{eq:gap}
  \gap X(S_n, Z_n) \le \gap X(S_n, S_{(n-2,2)}, Z_n) = 1 
  \qquad n\ge 3 \,.
\end{equation}
Our main result confirms what in \cite{GSW} was indicated as
``empirical evidence''.

\begin{theorem}\label{thm:main}
Let $Z_n := \{r_1, \ldots, r_n\}$ be the set of all initial reversals
defined in \eqref{eq:rk}.
Then, for any $n\ge 3$, we have $\gap X(S_n, Z_n)=1$.
\end{theorem}

\noindent
Our approach is based on the connection between
the Laplacian of a Cayley graph for a finite group $G$
and the irreducible representations of $G$.
A similar approach has allowed a detailed study
of the spectrum of Cayley graphs on $S_n$ generated by 
a set of transpositions $Z$,
when $Z$, interpreted as as the edge set of a graph
with $n$ vertices, yields a complete graph \cite{DiSh}, 
or a complete multipartite graph \cite{Ces1}.

\secskip
\section{Cayley graphs and irreducible representations} 
\label{sec:ind}

\noindent
In this section we introduce our notation and briefly recall
a few basic facts about the eigenvalues of (weighted) Cayley graphs
and the irreducible representations of a finite group.
Details can be found for instance in \cite{Lub2}.
Let $\rY$ be a representation of a finite group $G$. $\rY$ extends
to a representation of the complex \textit{group algebra} $\bC G$ by letting
$
  \rY(w) := 
  \sum_{g\in G} w_g \, \rY(g) 
$,
where $w = \sum_{g\in G} w_g \, g$ is an element in $\bC G$.
$\Irr(G)$ stands for the set of all equivalence classes
of irreducible representations
of $G$. If $[\a] \in \Irr(G)$ we denote with $\rT^\a$
a specific choice in the class $[\a]$.
By Maschke's complete reducibility theorem,
any representation $\rY$ is equivalent to a direct sum
\begin{equation}\label{eq:Y}
  \rY \cong \bigoplus_{[\a]\in \Irr(G)} y_\a \, \rT^\a \,,
\end{equation}
where $y_\a$ are suitable nonnegative integers.
By consequence the spectrum of $\rY(w)$
is just the union of the spectra of those $\rT^\a(w)$ for which $y_\a$
is nonzero.%
\footnote{if one is interested in multiplicities of the 
eigenvalues, spectra 
must be treated as multisets.}
We define the set of all \textit{symmetric} elements in $\bC G$ and
the set of all \textit{symmetric nonnegative} elements as
\begin{align*}
  \symCG &:= \{ w \in \bC G : w_g = \ol w_{g^{-1}} 
  \text{ for all $g\in G$} \} \\
  (\symCG)_+ &:= \{ w \in \symCG : w_g \in \bR,\ w_g \ge 0
  \text{ for all $g\in G$} \} \,.
\end{align*}
If $\rY$ is a \textit{unitary} representation and $w$ is symmetric,
then $\rY(w)$ is a Hermitian matrix. Since every representation
of a finite group is equivalent to a unitary representation,
the eigenvalues of $\rY(w)$ are real for \textit{any} representation $\rY$
and for any $w\in \symCG$. We denote with $\lmax( \rY(w))$ the
largest eigenvalue of $\rY(w)$.
A trivial upper bound on this quantity is found by
assuming $\rY$ unitary
\begin{equation}\label{eq:triv}
  \lmax(\rY(w)) \le \|\rY(w)\| \le \sum_{g\in G} |w_g| \, \|\rY(g)\| = 
  \sum_{g\in G} |w_g|  =: |w| \,,
\end{equation}
where $\|A\|$ stands for the $\ell^2$ operator norm of the matrix $A$
and $|w|$ for the $\ell^1$ norm on $\bC G$.

\medno
If $Z$ is a generating set for $G$ we can define an element of the group 
algebra $\bC G$, which we denote by $\hZ$, given by
\begin{equation*}
  \txt \hZ := \sum_{z\in Z} z \,.
\end{equation*}
In the following we will consider  \textit{symmetric} generating
sets $Z$, that is such that $z\in Z$ iff $z^{-1}\in Z$.
In this case $\hZ$ is an element of $\symCGp$.
Conversely if $w =\sum_{g\in G} w_g g$ is a symmetric nonnegative element
in $\bC G$, we can define the
(undirected) \textit{weighted Cayley graphs} $X(G,w)$,
where $w_g$
represents the weight associated to each edge $(h, gh)$, $h\in G$.
The adjacency matrix and the Laplacian of $X(G,w)$
are closely related to the \textit{left regular} representation
of $G$. If we denote with $\rR$ such a representation, it
follows from the definitions that
\begin{equation}\label{eq:azi1}
  A_{X(G,w)}  = |w|\, \id - \D_{X(G,w)} = \rR(w) 
  \qquad w \in \symCGp\,.
\end{equation}

\medno
Consider now the case in which $G$ is the permutation group $S_n$.
There is a one-to-one correspondence between
$\Irr(S_n)$ and the set of all partitions of $n$.
A \textit{partition} of $n$ is a nonincreasing 
sequence of positive integers 
$\a=(\a_1, \a_2, \ldots, \a_r)$ such that $\sum_{i=1}^r \a_i = n$.
We write $\a\partit n$ if $\a$ is a partition of $n$.
We denote with $[\a]$ the class of irreducible
representations of $S_n$ corresponding to the partition $\a$.
For simplicity we write $[\a_1, \ldots, \a_r]$ 
instead of $[(\a_1, \ldots, \a_r)]$.

Since all irreducible representations appear in the decomposition
of the left regular representation, it follows from \eqref{eq:azi1} that
if we let
\begin{equation}\label{eq:psi}
  \psi([\a], w) := |w| - \lmax( \rT^\a(w)) 
  \qquad \a\partit n \,,
\end{equation}
then the spectral gap of $X(S_n,w)$ is given by
\begin{equation}\label{eq:cay}
  \gap X(S_n,w)  =
  \min_{\a\partit n: \, \a\ne(n)}
  \psi( [\a] , w) \,,
\end{equation}
where $\a=(n)$ is the one--dimensional identity representation
which yields one eigenvalue equal to $|w|$ in $\rR(w)$.

\medno
We conclude this section with a remark concerning a connection between
results like Theorem \ref{thm:main} 
and  \textit{Aldous's conjecture} \cite{Ald} asserting 
that the random walk and the 
interchange process have the same spectral
gap on any finite graph.
In order to explain this connection
we introduce a property, which we call property (A) 
which is an attribute of certain elements of the group algebra:
given $w \in (\Smm \bC S_n)_+$, we say that property (A) holds
for $w$ if one of the following two equivalent statements is satisfied
\begin{enumerate}
\item[(A1)]
If $\a\partit n$ and $\a \ne (n)$, then
$\lmax( \rT^\a(w) ) \le \lmax( \rT^{(n-1,1)}(w) )$
\item[(A2)]
$ \gap X(S_n,w) = \psi(\, [n-1,1], w)$.
\end{enumerate}
The two statements are equivalent in virtue of \eqref{eq:cay}.
Aldous's conjecture, originally formulated in the
framework of continuous time Markov chains,
is equivalent (see \cite{Ces1}) to the assertion that:
\textit{if $w=\sum_{\ell} t_\ell$ is a sum of transpositions 
$t_\ell=(i_\ell \, j_\ell) \in S_n$,
then $w$ has property (A)}.
A stronger ``weighted graphs'' version of this conjecture can be formulated
in which $w=\sum_{\ell} w_\ell\, t_\ell$ is allowed to be a 
\textit{linear combination} of transpositions with nonnegative coefficients.
A weaker version of this statement, namely for bipartite graphs,
was also conjectured in \cite{Fri2}.
Several papers have appeared with proofs of Aldous's conjecture
for some particular classes of graphs, and recently a
beautiful general proof has been found by Caputo, Liggett and
Richthammer \cite{CaLiRi} (see also this paper for references to
previous work).
Going back to our problem of finding
the spectral gap of the Cayley graph $X(S_n, \hZ_n)$, where
$Z_n$ is the set of initial reversals, we observe that
Proposition 4.1 in \cite{GSW} implies that $\psi([n-1,1], \hZ_n) = 1$,
hence Theorem \ref{thm:main}
is equivalent to the assertion that $\hZ_n$ has property (A).

\secskip
\section{Proof of Theorem \ref{thm:main}} 

\noindent
We start with a general lower bound on the spectral gap
of a weighted Cayley graph $X(S_n,w_n)$ which makes use
of the \textit{branching rule} \cite[Section 2.8]{Sag}
for the decomposition of the restriction of an irreducible
representation $[\a]$ of $S_n$ to the subgroup $S_{n-1}$.
This rule states that
\begin{align*}
  [\a] \resR^{S_n}_{S_{n-1}} = \bigoplus_{\b \in \a^-} \, [\b]
  \qquad \a\partit n
\end{align*}
where, if $\a=(\a_1, \ldots,\a_r)$, 
 $\a^-$ is defined as the collection of all sequences
of the form
\begin{equation*}
  (\a_1, \ldots, \a_{i-1}, \a_i -1, \a_{i+1}, \ldots, \a_r)
\end{equation*}
\textit{which are partitions of} $n-1$. 
For example
\begin{align*}
  [6,5,5,3,1]\resR^{S_{20}}_{S_{19}} = 
  [5,5,5,3,1]
  \oplus [6,5,4,3,1]
  \oplus [6,5,5,2,1]
  \oplus [6,5,5,3] \,.
\end{align*}

\medno
We have then the following lower bound on the spectral
gap of $X(S_n,w_n)$.

\begin{lemma}\label{thm:dec}
Let $z_k \in (\Sym \bC S_k)_+$ for $k=1,2, \ldots$, and let 
$w_n := \sum_{k=1}^n z_k$. Then
\begin{equation}\label{eq:dec}
  \gap X(S_n, w_n) \ge 
  \min_{k = 2, \ldots, n} 
  \psi(\, [k-1,1],\, w_k ) 
  \,.
\end{equation}
\end{lemma}

\begin{remark}\label{thm:deco}
Given $w_n = \sum_{\pi\in S_n} w_{n,\pi}\, \pi \in (\Sym \bC S_n)_+$,
it is always possible to write $w_n$
as a sum of $z_k$ such that Lemma \ref{thm:dec} applies. For instance
one can define
\begin{equation*}
  z_k = \sum_{\pi \in S_k\setm S_{k-1}} w_{n,\pi}\, \pi \,,
\end{equation*}
even though it is not clear that this choice gives the optimal lower bound
in \eqref{eq:dec}.
\end{remark}

\begin{remark}\label{thm:rem}
Consider the case in which 
$w_n=\sum_\ell w_{n,\ell} \, t_\ell$ is a linear 
combination of transpositions $t_\ell \in S_n$ with $w_{n,\ell} \ge 0$,
and define the graph $\cG_{w_n}$ with vertex set $\{1,\ldots,n\}$
and edge set given by $\supp w_n = \{ t_\ell : w_{n,\ell} >0 \}$, in which
each transposition $t_\ell =(ij)$ is identified with 
the corresponding edge $\{i,j\}$.
In the case of transpositions Lemma \ref{thm:dec}
is equivalent to Lemma 2 in \cite{Mor} and it
was more or less implicit already in \cite{HaJu},
where it was used to prove Aldous's conjecture for trees,
meaning for all $w$ such that $\cG_w$ is a tree.
Using this approach, Aldous's conjecture
has been
proved independently in \cite{Mor} and \cite{StCo} for hypercubes 
\textit{asymptotically}, \ie in the limit when the side length
of the cube tends to infinity.
While the proof of Lemma 2 in \cite{Mor} 
(or the equivalent statement in \cite{HaJu}) is not hard,
it is nevertheless interesting to realize that our general 
formulation of this result
is a direct consequence of very general algebraic identities
(equality \eqref{eq:cay} and the branching rule).
\end{remark}

\noindent
\textit{Proof of Lemma \ref{thm:dec}}.
If $A$ and $B$ are two Hermitian $n\times n$ matrices we have
\begin{equation*}
  \lmax(A+B) = \max_{x\in \bC^n:\, \|x\|=1} \<(A+B) x, x\>
  \le \lmax(A) + \lmax(B) \,,
\end{equation*}
where $\|x\|$ is the Euclidean norm.
Using this fact and the trivial bound \eqref{eq:triv}, we find
\begin{align*}
    \lmax (\rT^\a(w_n)) &\le 
  \lmax (\rT^\a( w_{n-1} )) 
  + \lmax (\rT^\a(z_n)) \\
  &\le
  \lmax (\rT^\a( w_{n-1} )) 
  + |z_n|
  \,.
\end{align*}
Since $w_{n-1}\in (\Sym \bC S_{n-1})_+$, we can write
$  \rT^\a(w_{n-1}) \cong \bigoplus_{\b \in \a^-} \, \rT^\b(w_{n-1})$,
thus
\begin{equation*}
  \lmax (\rT^\a(w_{n-1})) = \max_{\b \in \a^-} \, \lmax (\rT^\b(w_{n-1})) \,.
\end{equation*}
It follows from the branching rule that 
if $\a\ne(n)$ and $\a\ne(n-1,1)$, then 
the trivial partition $(n-1)$ is not contained in $\a^-$.
By consequence
\begin{align}\label{eq:n-1}
  \max_{\substack{\a\partit n:\\ \a\ne (n),\, \a\ne(n-1,1)}}
  \lmax(\rT^\a( w_n))
  \le
  \max_{\substack{\b\partit n-1:\\ \b\ne(n-1)} }
  \lmax(\rT^\b(w_{n-1})) + |z_n| \,.
\end{align}
Since $z_k$ and $w_n$ have nonnegative components, we get
\begin{equation}\label{eq:wz}
  |w_n| := \sum_{\pi\in S_n} w_{n,\pi} = 
  \sum_{\pi\in S_n} (w_{n-1,\pi} + z_{n,\pi}) =
  |w_{n-1}| + |z_n| \,.
\end{equation}
From \eqref{eq:cay}, \eqref{eq:n-1} and 
\eqref{eq:wz} we obtain 
\begin{equation}\begin{split}
  &\gap X(S_n, w_n) = 
  \min_{\a\partit n: \, \a\ne(n)}
  \bigl[ \, |w_n| - \lmax( \rT^\a(w_n)) \, \bigr] \\
  \quad &\ge
  \min\Bigl\{ 
  \min_{\substack{\b\partit n-1:\\ \b\ne(n-1)} }
  \bigl[\, |w_{n-1}| - \lmax( \rT^\b (w_{n-1})) \, \bigr] \,,\,
  \psi([n-1,1],w_n) \Bigr\} \,.
\end{split}\end{equation}
Hence we get the recursive inequality
\begin{equation}\begin{split}\label{eq:rec}
  \gap X(S_n, w_n) \ge 
  \min \bigl\{ 
  \gap X(S_{n-1},w_{n-1}) \, ,\,
  \psi(\, [n-1,1],\, w_n )
  \bigl\} \,.
\end{split}\end{equation}
When $n=2$ we have $S_2=\{1,(12)\}$ and 
$w_2 = w_{2,1} \cdot 1 + w_{2,(12)} \cdot (12)$.
A trivial computation
yields 
\begin{equation}\label{eq:n=2}
  \gap X(S_2, w_2) =   \psi(\, [1,1],\, w_2 ) = 2 w_{2,(12)}
\end{equation}
which, together with \eqref{eq:rec}, implies the Lemma.
\qed

\medno
\textit{Proof of Theorem \ref{thm:main}}.
Let $r_k$ be the permutation which reverse the order of the
first $k$ positive integers
\begin{equation*}
  r_k : \{1,2,\ldots, n\} \too 
  \{ k,\, k-1,\, \ldots,\, 2,\,1,\,  k+1, \ldots, n\} \,,
\end{equation*}
and let $w_n := \sum_{k=1}^n r_k$.
Let also $\rD_n$ be the $n$--dimensional defining representation of $S_n$,
which can be written as
\begin{equation*}
  [\rD_n] = [n] \oplus [n-1,1] \,.
\end{equation*}
The eigenvalues and eigenvectors of the matrix $\rD_n(w_n)$
are listed in \cite[Proposition 4.1]{GSW}. The two largest eigenvalues
are $n$ and $n-1$, both simple. Since the eigenvalue $n$
clearly corresponds to the identity representation $[n]$
contained in $\rD_n$, we have, \textit{for each} $n \ge 3$, 
\begin{equation*}
  \lmax( \rT^{(n-1,1)}(w_n)) = n-1\,.
\end{equation*}
By consequence,
\begin{align*}
  \psi(\, [n-1,1],\, w_n ) = 1 
  \qquad n\ge 3 \,.
\end{align*}
From Lemma \ref{thm:dec}, \eqref{eq:gap} 
and from \eqref{eq:n=2} which in this case says
$\psi(\, [1,1],\, w_2 ) = 2$,
it follows that $\gap(S_n,w_n) = 1$.
\qed

%
%
%
%

\begin{thebibliography}{10}

\bibitem{Ald}
David Aldous, \emph{\upshape
  www.stat.berkeley.edu/\~{}aldous/research/op/sgap.html}.

\bibitem{CaLiRi}
Pietro Caputo, Thomas~M. Liggett, and Thomas Richthammer, \emph{A recursive
  proof of {A}ldous' spectral gap conjecture}, \upshape arXiv:0906.1238v3
  (2009).

\bibitem{Ces1}
Filippo Cesi, \emph{On the eigenvalues of {C}ayley graphs on the symmetric
  group generated by a complete multipartite set of transpositions}, \upshape
  arxiv:0902.0727v1 (2009).

\bibitem{DiSh}
Persi Diaconis and Mehrdad Shahshahani, \emph{Generating a random permutation
  with random transpositions}, Z. Wahrsch. Verw. Gebiete \textbf{57} (1981),
  no.~2, 159--179.

\bibitem{Fri2}
Joel Friedman, \emph{On {C}ayley graphs on the symmetric group generated by
  transpositions}, Combinatorica \textbf{20} (2000), no.~4, 505--519.

\bibitem{GSW}
Paul~E. Gunnells, Richard~A. Scott, and Byron~L. Walden, \emph{On certain
  integral {S}chreier graphs of the symmetric group}, Electron. J. Combin.
  \textbf{14} (2007), no.~1, Research Paper 43, 26 pp. (electronic).

\bibitem{HaJu}
Shirin Handjani and Douglas Jungreis, \emph{Rate of convergence for shuffling
  cards by transpositions}, J. Theoret. Probab. \textbf{9} (1996), no.~4,
  983--993.

\bibitem{Lub2}
Alexander Lubotzky, \emph{Cayley graphs: eigenvalues, expanders and random
  walks}, Surveys in combinatorics, 1995 ({S}tirling), London Math. Soc.
  Lecture Note Ser., vol. 218, Cambridge Univ. Press, Cambridge, 1995,
  pp.~155--189.

\bibitem{Mor}
Ben Morris, \emph{Spectral gap for the interchange process in a box}, Electron.
  Commun. Probab. \textbf{13} (2008), 311--318.

\bibitem{Nas}
D.~Nash, \emph{Cayley graphs of symmetric groups generated by reversals}, Pi Mu
  Epsilon Journal \textbf{12} (2005), 143--147.

\bibitem{Sag}
Bruce~E. Sagan, \emph{The symmetric group: Representations, combinatorial
  algorithms, and symmetric functions}, second ed., Graduate Texts in
  Mathematics, vol. 203, Springer-Verlag, New York, 2001.

\bibitem{StCo}
Shannon Starr and Matt Conomos, \emph{Asymptotics of the spectral gap for the
  interchange process on large hypercubes}, \upshape arxiv:0802.1368v2 (2008).

\end{thebibliography}

\providecommand{\bysame}{\leavevmode\hbox to3em{\hrulefill}\thinspace}
\providecommand{\MR}{\relax\ifhmode\unskip\space\fi MR }
\providecommand{\MRhref}[2]{%
  \href{http://www.ams.org/mathscinet-getitem?mr=#1}{#2}
}
\providecommand{\href}[2]{#2}

\end{document}